\documentclass[10pt, reqno]{amsart}
\usepackage{amssymb,latexsym}
\usepackage{eucal}

\newtheorem*{proposition*}{}

\newtheorem*{Lem17.3}{Lemma 17.3}
\newtheorem*{Lem18.3}{Lemma 18.3}
\newtheorem*{Lem18.5}{Lemma 18.5}
\newtheorem*{Lem19.0}{Lemma 19.0}
\newtheorem*{Lem19.1}{Lemma 19.1}
\newtheorem*{Lem19.1.3}{Lemma 19.1.3}

\newtheorem*{Thm}{Theorem}

\DeclareMathOperator{\SQ}{\textup{SQ}}

\newcommand{\ph}{{\varphi}}
\newcommand{\D}{\Delta}

\newcommand{\A}{\mathcal{A}}
\newcommand{\B}{\mathcal{B}}

\newcommand{\HH}{\mathcal{H}}
\newcommand{\K}{\mathcal{K}}
\newcommand{\R}{\mathcal{R}}
\newcommand{\Ss}{\mathcal{S}}

\newcommand{\X}{\mathcal{X}}

\newcommand{\p}{\partial}
\newcommand{\ra}{\rangle}
\newcommand{\la}{\langle}

\begin{document}

\title[On subgroups of free Burnside groups]
{On subgroups of free Burnside groups of large odd exponent}
\author{S.V. Ivanov}
\address{Department of Mathematics\\
University of Illinois \\
Urbana,  IL 61801}
\email{ivanov@math.uiuc.edu}
\thanks{Supported in part by NSF grant DMS 00-99612}
\subjclass[2000]{Primary  20E07, 20F05, 20F50}

\begin{abstract}
We prove that every noncyclic subgroup of a free $m$-generator Burnside group
$B(m,n)$ of odd exponent $n \gg 1$ contains a subgroup $H$ isomorphic to a
free Burnside group $B(\infty,n)$ of exponent $n$ and countably infinite rank
such that for every normal subgroup $K$ of $H$ the normal closure $\langle K
\rangle^{B(m,n)}$ of $K$ in $B(m,n)$ meets $H$ in $K$. This implies that
every noncyclic subgroup of  $B(m,n)$ is SQ-universal in the class of groups
of exponent $n$.
\end{abstract}
\maketitle

A group $G$ is called {\em $\SQ$-universal} if every countable group is isomorphic
to a subgroup of a quotient of $G$. One of classical embedding theorems proved
by Higman, B. Neumann, H. Neumann in \cite{HNN49} states that every
countable group $G$ embeds in a 2-generator group or, equivalently, a free group
$F_2$ of rank 2 is $\SQ$-universal. Recall that the proof of this
theorem makes use of the following natural definition.
A subgroup $H$ of a group $G$ is  called a {\em $Q$-subgroup}
if for every normal subgroup $K$ of $H$ the normal closure
$\langle K \rangle^{G}$ of $K$ in $G$ meets $H$ in $K$, i.e.,
$\langle K \rangle^{G} \cap H = K$.
For example, factors $G_1$, $G_2$ of the free product $G_1 * G_2$
or the direct product $G_1 \times G_2$ are
$Q$-subgroups of $G_1 * G_2$ or $G_1 \times G_2$, respectively.
In particular, a free group $F_m$ of
rank $m>1$, where $m =\infty$ means countably infinite rank, contains
a $Q$-subgroup isomorphic to $F_k$ for every $k \le m$.
On the other hand, it is proved in \cite{HNN49} that the subgroup
$
\langle    a^{-1} b^{-1} a b^{-i} a b^{-1} a^{-1} b^{i}  a^{-1} b
a b^{-i}  a b  a^{-1} b^{i}  \; | \; i =1,2, \dots  \rangle
$
of $F_2 = F_2(a,b)$ is a $Q$-subgroup of $F_2$ isomorphic to $F_\infty$ and
freely generated by indicated elements.
In \cite{NN59} B. Neumann and H. Neumann found simpler generators and proved
that $\langle   [ b^{-2i+1} a b^{2i-1}, a]  \; | \; i =1,2, \dots  \rangle$,
where $[x,y] = xyx^{-1}y^{-1}$ is the commutator of $x$ and $y$, is a
$Q$-subgroup of $F_2$ isomorphic to $F_\infty$ and freely generated by
indicated elements.
It is obvious that the property of being a $Q$-subgroup is transitive.
Therefore, a group $G$ contains a
$Q$-subgroup isomorphic to $F_\infty$ if and only if $G$ contains a
$Q$-subgroup isomorphic to $F_m$, where $m \ge 2$.

Ol'shanskii \cite{O95} proved that any nonelementary subgroup of a hyperbolic
group $G$ (in particular, $G = F_m$) contains a $Q$-subgroup
isomorphic to $F_2$. In particular, if $G$ is a nonelementary hyperbolic  group
then $G$ is $\SQ$-universal.

It follows from an embedding theorem of Obraztsov (see Theorem 35.1 in \cite{O89})
that any countable group of odd exponent  $n \gg 1$ embeds in a 2-generator group
of exponent $n$ and so a free  2-generator Burnside group $B(2,n) = F_2 /F_2^n$
of exponent $n$ is $\SQ$-universal in the class of groups of exponent $n$.
Interestingly, the proof of this theorem has nothing to do with
free $Q$-subgroups of  the Burnside group $B(2,n)$
and does not imply the existence of such subgroups in $B(2,n)$.

Ol'shanskii and Sapir proved in \cite{OS02} (among many other things)
that for odd $n \gg 1$ the group
$B(m,n)$ with some $m = m(n)$ does contain  $Q$-subgroups isomorphic to
$B(\infty,n) = F_\infty / F_\infty^n$.
Sonkin \cite{S02} further refined corresponding arguments of \cite{OS02}
to show that for odd  $n \gg 1$ the group $B(2,n)$ contains
a $Q$-subgroup isomorphic to $B(\infty,n)$ which also implies that
$B(2,n)$ is $\SQ$-universal in the class of groups of exponent $n$.

Recall that an embedding $B(\infty,n) \to B(2,n)$ for odd $n \ge 665$,
without the property of being $Q$-subgroup, was first proved by
Shirvanian \cite{Sh76}.
Atabekian  \cite{A86}, \cite{A87} showed
for odd $n \gg 1$ (e.g., $n > 10^{78}$) that every noncyclic subgroup
of $B(m,n)$ contains a subgroup isomorphic  to $B(2,n)$ (and so,
by Shirvanian's result, contains a subgroup isomorphic  to $B(\infty,n)$).
A short proof of this Atabekian's theorem due to the author was incorporated
in \cite{O89} (Theorem 39.1). It turns out that the same idea of "fake" letters
and using relations of Tarski monsters yields not only embeddings but also
embeddings as  $Q$-subgroups and significantly shortens corresponding
arguments of \cite{OS02}, \cite{S02}. The aim of this note is to elaborate
on this idea and to strengthen Atabekian's theorem as follows.

\begin{Thm} Let $n$ be odd, $n \gg 1$ (e.g., $n >10^{78}$),
and $B(m,n)$ be a free $m$-generator Burnside group of exponent $n$.
Then every noncyclic subgroup of $B(m,n)$ contains a $Q$-subgroup of $B(m,n)$
isomorphic to $B(\infty,n)$. In particular, every noncyclic subgroup of
$B(m,n)$ is $\SQ$-universal in the class of groups of exponent $n$.
\end{Thm}

{\em Proof of Theorem.}
To be consistent with the notation of \cite{O89}, rename
the exponent $n$ by $n_0$. Consider an alphabet $\A = \{ a_1, \dots, a_m\}$ with
$m \ge 2$. Let $G(\infty)$ be a presentation for the free
Burnside group $B(\A,n_0)$ of exponent $n_0$ in the alphabet $\A$
constructed  as in Sect. 18.1 \cite{O89} and $\HH$ be a noncyclic subgroup
of $B(\A,n_0)$. Conjugating if necessary, by Lemma 39.1 \cite{O89},
we can suppose that there are words $F, T \in \HH$ such that $F$ is a period
of some rank $|F|$ (with respect to the presentation $G(\infty)$ of
$B(\A,n_0)$),  $|T| < 3|F|$ and $FT \neq TF$ in $B(\A,n_0)$.

Consider a presentation
\begin{equation}
\K = \langle b_1, b_2 \;  \|  \;  R=1, R \in \bar \R_0 \rangle
\end{equation}
for a 2-generator
group $\K$ of exponent $n_0$ (perhaps, $\K$ is trivial).

Set $\bar \A = \A \cup \{ b_1, b_2 \}$ and define
$
\bar G(0) = \langle \bar \A \; \| \; R=1, R \in \bar  \R_0 \rangle .
$
Clearly, $\bar G(0)$ is the free product of the free group $G(0) = F(\A)$
in $\A$ and $\K$. If $W$ is a word in $\bar \A^{\pm 1} = \A \cup \A^{-1}$
then its {\em length} $|W| = |W|_\A$ is defined to be the number of letters
$a_k^{\pm 1}$, $a_k \in \A$, in $W$.
In particular, $| b_1 |= | b_2|= 0$. Using this new length, we construct
groups
$
\bar G(i) =  \langle \bar \A \; \| \; R=1, R \in \bar \R_i \rangle \
$
by induction on $i\ge 1$ exactly as in Sect. 39.1 \cite{O89}, that is, the set
$\bar \Ss_i = \bar  \R_i \setminus \bar \R_{i-1}$ of defining words of rank $i$
consists of all relators
of the first type $A^{n_0}$, $A \in \bar \X_i$, if $i < |F|$.
As before, we observe
that the set $\bar \X_{|F|}$ of periods of rank $i = |F|$ can be chosen
so that $F \in \bar \X_{ |F|}$. For $i = |F|$ the set
$\bar \Ss_i = \bar  \R_i \setminus \bar \R_{i-1}$ consists of all relators
of the first type $A^{n_0}$, $A \in \bar \X_i$, and two relators of the second type
which are words of the form
\begin{equation}
b_1 F^{n} T F^{n+2} \dots T F^{n+2h -2} , \quad
b_2 F^{n+1} T F^{n+3} \dots T F^{n+2h -1} .
\end{equation}

For $i > |F|$ the set $\bar \Ss_i = \bar  \R_i \setminus  \bar  \R_{i-1}$
again consists of all relators $A^{n_0}$
of the first type only, $A \in \bar \X_i$. Thus, the groups
$$
\bar G(i)=  \langle \bar \A \; \| \; R=1, R \in \bar \R_i \rangle , \qquad
\bar G(\infty) = \langle \bar \A \; \| \; R=1, R \in
\cup_{j=0}^\infty \bar \R_j \rangle
$$
are constructed.

Consider a modification of condition $R$ (Sect. 25.2 \cite{O89}),
that will be called  {\em condition $R'$}, in which property $R4$ is replaced by
the following.
\begin{enumerate}
\item[$R4'$]  The words $T_k$ are not contained in the subgroup
$\la A \ra$ of the group  $\bar G(i-1)$, $i \ge 1$,  except for the case
when $k=1$, $|T_1| = 0$ and the integers $n_1, n_k$ have the same sign.
\end{enumerate}

Let $\bar \D$ be a diagram over the graded presentation $\bar G(i)$, $i \ge 0$.
According to the new definition of the word length, we define
the length $|p|$ of a path $p$ so that $|p| = |\ph(p)|$. In particular, if $e$
is an edge of $\bar \D$  with $\ph(e) = b_k^{\pm 1}$, $k =1,2$, then $|e| =0$.
Hence, such an edge $e$ is regarded to be a 0-edge of $\bar \D$  of {\em type 2}.
Recall that if $\ph(e) =1$ then $e$ is called in  \cite{O89}
a 0-edge (we will specify that such $e$ is a 0-edge of {\em type 1}).
All faces labelled by relators of $\bar G(0)$ are also called 0-{\em faces}
(or faces of rank 0) of $\bar \D$.  A 0-face $\Pi$ of $\bar \D$ has {\em type 1}
if it is a 0-face in the sense of \cite{O89}. Otherwise, i.e., when $\p \Pi$
has a nontrivial label $R \in \bar \R_0^{\pm 1}$, a 0-face $\Pi$ has {\em type 2}.

Note that the new definition of length and the existence of 0-edges of type 2
imply a number of straightforward corrections  in analogues of definitions
and lemmas of Sects. 18, 25 of
\cite{O89} on  group presentations with condition  $R'$ (these corrections
are quite analogous to what is done in similar situations of papers
\cite{I02a}, \cite{I02b}).  For example, in the definition of a simple in rank $i$
word $A$ (Sect. 18.1 \cite{O89}) it is in addition required that $|A| > 0$.
Lemma 25.1 now claims that every reduced diagram
$\bar \D$ on a sphere or torus  has rank 0. Corollary 25.1 is stated
for $\bar \R_i \setminus \bar \R_0$  and Corollary 25.2 is now missing.
In Lemma 25.2 we additionally allow that $X$ is conjugate to a word of length 0.
Lemmas 25.12--25.15 are left out.

Repeating the proof of Lemma 27.2 \cite{O89} (and increasing the number
of short sections in Lemma 27.1 \cite{O89} from 2 to 3), we can show that
presentations   $\bar G(i)$, $\bar G(\infty)$ satisfy condition $R'$.
Furthermore, it is straightforward to check that proofs of Lemmas
26.1--26.5 for a graded presentation with  condition $R'$ are retained
(minor changes in arguments of proofs of Lemmas 26.1--26.2 caused by the
possibility $|T_1| = 0$ are evident). Thus, by Lemma 26.5 \cite{O89}, any reduced
diagram over $\bar G(i)$ (or $\bar G(\infty)$) is a  $B$-map.

By definitions and the analogue of Lemma 25.2 \cite{O89},
the group $\bar G(\infty)$  has exponent $n_0$.
Suppose $U$ is a word in $\{ b_1^{\pm 1}, b_2^{\pm 1} \}$ and $U =1$ in
the group $\bar G(\infty)$. Let $\bar \D$
be a reduced diagram over $\bar G(\infty)$   with $\ph(\p \bar \D) \equiv U$.
Since $|\p \bar \D | = 0$, it follows from Theorem 22.4 \cite{O89} that
$r(\bar \D) = 0$, hence $U =1$ in the group $\K$ given by (1).
This means that  $\K$ naturally embeds in $\bar G(\infty)$.

Let $V_1 = (F^n T F^{n+2} \dots T F^{n+2h -2} )^{-1}$,
$V_2 = (F^{n+1} T F^{n+3} \dots T F^{n+2h -1})^{-1}$.
Observe that, in view of relators (2),
the group $\bar G(\infty)$ is naturally isomorphic to the
quotient
$
B_\K(\A, n_0) = \la   B(\A, n_0)   \; \| \; R(V_1, V_2)=1, \;
R(b_1, b_2) \in  \bar \R_0  \ra
$
of $B(\A, n_0)$. Hence, the subgroup $\la V_1, V_2 \ra$ of $B_\K(\A, n_0)$
is isomorphic to the group $\K$ given by (1) under the map
$V_1 \to b_1$, $V_2 \to b_2$. Since $\K$ is an arbitrary 2-generator
group of exponent $n_0$, it follows that  $\la V_1, V_2 \ra$ is a
$Q$-subgroup of  $B(\A, n_0)$ isomorphic to $B(2, n_0)$.
\smallskip

Now we will show that   $B(\infty, n_0)$ embeds in  $B(2, n_0)$ as a $Q$-subgroup.
To do this we will repeat the above arguments with some changes.
Now we let $\A = \{ a_1, a_2 \}$, i.e., $m =2$, and
$\B = \{ b_1, b_2, \dots \}$ be a countably infinite alphabet. Let
\begin{equation}
\K = \la \B \; \| \; R=1, R \in \bar \R_0 \ra
\end{equation}
be a presentation of a finite or countable group  of exponent $n_0$,
$\bar \A = \A \cup \B$, and
$
\bar G(0) = \la \bar \A \; \| \; R=1, R \in \bar \R_0 \ra .
$

As before, constructing groups $\bar G(i)$ by induction on $i \ge 1$,
we first define the set $\bar \X_i$ of periods of rank $i \ge 1$.
It is easy to show that each $\bar \X_i$, $i \ge 1$,  contains a word $A_i$ in
the alphabet $\{ a_1, a_2\}$ such that $A_i$ is not in the cyclic
subgroup $\la a_1 \ra$ of $\bar G(i-1)$.  Then for every
$F \in  \bar \X_i$ we define the relator $F^{n_0}$ and for the distinguished
period $A_i \in  \bar \X_i$ we introduce the second relator
$$
b_i A_i^{n} a_1 A_i^{n+2} \dots a_1 A_i^{n+ 2h -2} .
$$
These relators over all $F \in \bar \X_i$ form the set
$\bar \Ss_i = \bar \R_i \setminus \bar \R_{i-1}$. As above, we set
$$
\bar G(i) =   \la \bar \A \; \| \; R=1, R \in \bar \R_i \ra , \qquad
\bar G(\infty) =   \la \bar \A \; \| \;
R=1, R \in \cup_{j=0}^\infty \bar \R_j \ra
$$
and show that these presentations satisfy condition $R'$.
Similarly, we establish analogues of corresponding claims of Sects.
18, 25--27 \cite{O89}.

Suppose $U = U(\B)$ is a word in $\B^{\pm 1}$ and $U =1$ in
the group $\bar G(\infty)$.  Let $\bar \D$
be a reduced disk diagram over $\bar G(\infty)$  with $\ph(\p \bar \D) \equiv U$.
It  follows from Lemma 26.5, Theorem 22.4 \cite{O89}
and equality $|\p \bar \D | = 0$ that $r(\bar \D) = 0$. Hence $U =1$ in the group $\K$
given by (3) and so $\K$ naturally embeds in $\bar G(\infty)$.
As above, by definitions and Lemma 25.2 \cite{O89}, the group $\bar G(\infty)$
has exponent $n_0$ and we can see that $\bar G(\infty)$ is  naturally
isomorphic to the quotient
$$
B_\K(\A, n_0) = \la   B(\A, n_0)   \; \| \; R(V_1, V_2, \dots)=1, \;
R(b_1, b_2, \dots) \in \bar \R_0  \ra
$$
of $B(\A, n_0)$, where
$V_i = ( A_i^{n} a_1 A_i^{n+2} \dots a_1 A_i^{n+ 2h -2} )^{-1}$,
$i =1,2,\dots$.
Hence, the subgroup $\la V_1, V_2, \dots \ra$ of $B_\K(\A, n_0)$
is isomorphic to the group  $\K$ under the map
$V_i \to b_i$, $i =1,2, \dots$. Since $\K$ is an arbitrary
finite or countable group of
exponent $n_0$, it follows that  $\la V_1, V_2, \dots \ra$ is a
$Q$-subgroup of  $B(\A, n_0)= B(2, n_0)$ isomorphic to $B(\infty, n_0)$.

The specific estimate $n = n_0 > 10^{78}$ of Theorem can be obtained
by using lemmas and specific estimates of articles \cite{O82} and \cite{AI87}
(see also \cite{O85}) instead of those of \cite{O89}.
The proof of Theorem is complete. \qed
\smallskip

In conclusion, we remark that it is not difficult
to show that $B(\infty,n)$ embeds
in $B(2,n)$ for $n = 2^k \gg 1$ (see \cite{IO97}, \cite{I94})
but  it is not clear how to embed $B(\infty,n)$  in $B(2,n)$
as a $Q$-subgroup and it would be
interesting to do so. It would also be of interest to find out
whether  $B(\infty,n)$  embeds (as a $Q$-subgroup) in every nonlocally
finite subgroup of $B(m,n)$ for $n = 2^k \gg 1$.

\end{document}